\documentclass{amsart}
\usepackage{amsmath,amssymb,amscd,latexsym}
\usepackage[multiple]{footmisc}
\usepackage[all]{xy}

\newcounter{number}[section]

\newenvironment{nummer}{\refstepcounter{number}{\noindent\arabic{section}.\arabic{number}}}{}

\newcommand{\bn}{\noindent \begin{nummer} \rm}
\newcommand{\en}{\end{nummer}}

\newenvironment{ntheorem}{\noindent {\sc Theorem:} \it}{}
\newenvironment{nlemma}{\noindent {\sc Lemma:} \it}{}
\newenvironment{nprop}{\noindent {\sc Proposition:} \it}{}
\newenvironment{ndefn}{\noindent {\sc Definition:} \it}{}
\newenvironment{ncor}{\noindent {\sc Corollary:} \it}{}
\newenvironment{nconj}{\noindent {\sc Conjecture:} \it}{}

\newenvironment{nquestion}{\noindent {\sc Question:} }{}

\begin{document}

\title[Nuclear $\mathrm{C}^{*}$-algebras]{{\sc Structure of nuclear $\mathrm{C}^{*}$-algebras:\\From quasidiagonality to classification,\\and back again}}

\author{Wilhelm Winter}
\address{Mathematisches Institut\\
WWU M\"unster\\ Germany}

\email{wwinter@uni-muenster.de}

\date{November 30, 2017}
\thanks{Supported by the Deutsche Forschungsgemeinschaft through SFB 878.}

\begin{abstract}
I give an overview of recent developments in the structure and classification theory of separable, simple, nuclear C$^*$-algebras. I will in particular focus on the role of quasidiagonality and amenability for classification, and on the regularity conjecture and its interplay with internal and external approximation properties. 
\end{abstract}

\maketitle

\section*{Introduction}

\noindent
A C$^*$-algebra is a (complex) Banach-$^*$ algebra such that $\|x^*x\| = \|x\|^2$ for all elements $x$. Equivalently, C$^*$-algebras may be thought of as norm-closed $^*$-subalgebras of the bounded operators on Hilbert spaces. A von Neumann algebra is one that is even closed with respect to the weak operator topology. Examples of C$^*$-algebras include continuous functions on compact Hausdorff spaces, section algebras of vector bundles with matrix fibres, or suitable norm completions of group algebras. Group C$^*$-algebras come in different sizes; in particular there is a full one, which is universal with respect to all unitary representations of the group, and a reduced one, which is the norm completion of the left regular representation. Similar constructions can be associated with topological dynamical systems, via the crossed product construction.  

From the 1970s on it became clear that the notion of amenability for groups, with its many equivalent formulations, can be rephrased, in almost as many ways, for operator algebras as well. Some of these notions are more or less directly carried over from groups to group C$^*$-algebras or group von Neumann algebras. But then it often turns out that they make perfect sense at a much more abstract level --- and even there they remain closely related. Highlights are Choi--Effros' and Kirchberg's characterisation of nuclear C$^*$-algebras by the completely positive approximation property, and of course Connes' classification of injective $\mathrm{II}_1$ factors.  

Connes' theorem kicked off an avalanche of further developments in von Neumann algebras, but it also remained an inspiration for C$^*$-algebras. Elliott quite boldly conjectured that separable nuclear C$^*$-algebras should be classifiable by K-theoretic data. After some refinements and adjustments, and many years of hard work, we now understand the conjecture and its scope fairly well, at least for simple and unital C$^*$-algebras. Moreover, after a long detour, we now know that classification of simple nuclear C$^*$-algebras is not only philosophically, but also technically, surprisingly analogous to the classification of injective factors. In particular, the von Neumann algebraic properties of being (o) injective, (i) hyperfinite, (ii) McDuff, and (iii) having tracial comparability of projections, do have C$^*$-algebraic counterparts: (o) nuclearity, (i) finite (noncommutative) topological  dimension, (ii) $\mathcal{Z}$-stability, and (iii) comparison of positive elements. Here, nuclearity is characterised via the completely positive approximation property, topological dimension is either decomposition rank or nuclear dimension, and $\mathcal{Z}$-stability is tensorial absorption of the Jiang--Su algebra $\mathcal{Z}$, the smallest possible C$^*$-algebraic analogue of the hyperfinite $\mathrm{II}_1$ factor $\mathcal{R}$. Comparison of positive elements can be described as a regularity property (lack of perforation) of the Cuntz semigroup. 

A major difference between the C$^*$- and the von Neumann algebra situation is that for von Neumann factors injectivity implies the three other properties, whereas a simple C$^*$-algebra may be nuclear but fail to have finite nuclear dimension, be $\mathcal{Z}$-stable, or have strict comparison. However, at least conjecturally these three properties occur or fail simultaneously --- and for C$^*$-algebras with, say, not too complicated tracial state spaces, this is indeed a theorem. For now let us state a special case; we give a more comprehensive version later on. 

\noindent
{\sc Theorem A:} \emph{For a separable, simple, unital, nuclear $\mathrm{C}^*$-algebra $A \neq M_r(\mathbb{C})$ with at most one tracial state, the following are equivalent:
\begin{itemize}
\item[(\rm{i})] $A$ has finite nuclear dimension.
\item[(\rm{ii})] $A$ is $\mathcal{Z}$-stable, $A \cong A \otimes \mathcal{Z}$.
\item[(\rm{iii})] $A$ has strict comparison of positive elements. 
\end{itemize}
}

After the classification of injective factors was complete, it took around fifteen years to finish  the C$^*$-analogue of the type III case. (I say `finish', but this is only correct modulo the UCT problem. We will soon return to this little wrinkle.) So why did the type II analogue take almost forty years? For once, there is more information to keep track of: for type $\mathrm{II}_1$ factors the invariant is simply a point, whereas for C$^*$-algebras the invariant involves all possible ordered K-groups together with arbitrary Choquet simplices (arising as tracial state spaces). There is, however, also a deeper, and more mysterious  reason. This is related to both the universal coefficient theorem (UCT) and to quasidiagonality. The UCT problem asks whether all separable nuclear C$^*$-algebras are---in a very weak sense---homotopy equivalent to commutative ones. Conceptually this has a topological flavour, so from this perspective it is reasonable that the UCT becomes an issue in the C$^*$-algebraic (i.e., topological) setup, and not in the von Neumann algebraic (i.e., measure theoretic) situation. Nonetheless, I would like to understand this explanation at a more technical level --- but maybe this is asking for too much as long as the UCT problem is not resolved.  Quasidiagonality is an external approximation property; the quasidiagonality question (QDQ) asks whether all stably finite nuclear C$^*$-algebras admit a separating set of finite dimensional approximate representations; cf.\ \cite{BlaKir:MathAnn}. For von Neumann algebras the situation is more clear, since here it is a 2-norm (i.e., \emph{tracial}) version of quasidiagonality that matters.  

The connection between amenability and quasidiagonality was perhaps first drawn in \cite{Hadwin:JOT}, where Rosenberg observed that discrete groups with quasidiagonal reduced group C$^*$-algebras are amenable. The converse statement became known as Rosenberg's conjecture.

In \cite{TWW}, QDQ was answered for UCT C$^*$-algebras with faithful tracial states. 

\noindent
{\sc Theorem B:} \emph{Let $A$ be a nuclear $\mathrm{C}^*$-algebra with a faithful tracial state. Suppose $A$ satisfies the UCT (one could also say that $A$ is $\mathrm{KK}$-equivalent to a commutative $\mathrm{C}^*$-algebra). 
Then, $A$ is quasidiagonal.}

By work of Tu (and Higson--Kasparov), amenable group C$^*$-algebras do satisfy the UCT; since they also have a canonical faithful trace, this confirms Rosenberg's conjecture, and one arrives at a new characterisation of amenability.

\noindent
{\sc Corollary C:} \emph{For a discrete group $G$, the reduced group $\mathrm{C}^*$-algebra is quasidiagonal if and only if $G$ is amenable.}

Upon combining Theorems A and B with the work of Elliott, Gong, Lin, and Niu, we are now in a position to state the most general classification result that can currently be expected in the simple and unital case.\footnote{To be precise, we need more general versions of Theorems A and B here; I'll state these in Sections~{\ref{sectionTW}} and {\ref{section-qd}.}}\footnote{There are impressive results also in the non-unital and even in the non-simple situation, but I won't go into these here.}

\noindent
{\sc Theorem D:} \emph{The class of separable, simple, unital, nuclear, and $\mathcal{Z}$-stable $\mathrm{C}^*$-algebras satisfying the UCT is classified by $\mathrm{K}$-theoretic invariants.}

Separability will always be necessary for a classification result of this type, and nuclearity and $\mathcal{Z}$-stability are known to be essential --- and so, within its simple and unital scope, the theorem is complete modulo the UCT problem. For the time being this remains a sore point; however, one should note that in applications the C$^*$-algebras of interest  very often come with sufficient additional geometric structure, so that the UCT can be verified directly.

With the benefit of hindsight, one might divide the classification programme into three major challenges.  
These are linked at many levels---not least by the final result---but \emph{quasidiagonality} showcases these connections particularly beautifully.  This is the point of view I will take in these notes. Let us have a quick look at each of these challenges and very briefly sketch how quasidiagonality enters the game; we will see some more details in the main part of the paper.

\noindent
The first challenge: \emph{Understand nuclearity and the interplay with finite dimensional approximation properties.}\\ 
The main step was the completely positive approximation property as established by Choi--Effros and Kirchberg. This has been refined in various ways since then; in particular it has been used to model finite covering dimension in a noncommutative setting. While these are \emph{internal} approximations, quasidiagonality may be regarded as an \emph{external} approximation property. Understanding when quasidiagonality holds is a major task of the theory.

\noindent
The second challenge: \emph{Understand the $\mathrm{C}^*$-algebraic regularity properties} (i), (ii), (iii)  \emph{above and their interplay.} \\
This is about the regularity conjecture for nuclear C$^*$-algebras, now often referred to as Toms--Winter conjecture. I will state the conjecture in its full form, and describe what we know and what we don't know. There are two C$^*$-algebraic counterparts of hyperfiniteness in this context: finite decomposition rank and finite nuclear dimension. The first occurs only for finite C$^*$-algebras, the second in greater generality. It was open for some time what the difference between the two notions is, and we now know that (at least for simple C$^*$-algebras) the dividing line is marked by quasidiagonal traces.

\noindent
The third challenge: \emph{Implement the actual classification procedure.} \\
This is technical, and not easy to describe in short. As an illustration I will at least state a stable uniqueness result, Theorem~\ref{su-theorem}, which allows one to compare $^*$-homomorphisms up to unitary equivalence. This is particularly important for Elliott's intertwining argument; cf.\ \cite{Ell:intertwining-in-tuxedo}. I will also mention Kirchberg--Phillips classification,  Lin's TAF classification, and the recent spectacular work of Gong--Lin--Niu, together with the, indeed quite final, classification theorem that is now in place. 

I find the common history of quasidiagonality and classification most intriguing, since the two have met time and time again, and often very unexpectedly. First, Kirchberg used Voiculescu's result on quasidiagonality of suspensions to prove his famous $\mathcal{O}_2$-embedding theorem, which in turn led to Kirchberg--Phillips classification of purely infinite nuclear C$^*$-algebras. Next, Popa showed how to excise finite dimensional C$^*$-algebras in quasidiagonal ones with many projections. This was the starting point for Lin's TAF classification. (Later on, quasidiagonality also became crucial for the classification of specific types of examples, in particular for simple quotients of certain group C$^*$-algebras and for crossed products associated to free and minimal $\mathbb{Z}^d$-actions on compact and finite dimensional Hausdorff spaces.) Quasidiagonality (of all traces, to be precise) was also a crucial hypothesis for the classification result of Elliott--Gong--Lin--Niu. Theorem B above marked a surprising turn of events, as it, conversely, invoked a classification result to arrive at quasidiagonality.

The title of this note refers to the chronological development of matters, as outlined above. The main body of the paper is arranged thematically, in order to give a better overview of the individual aspects of the theory.

In Section~{\ref{section1}} we recall the notion of nuclearity, and various versions of the completely positive approximation property. Section~{\ref{section2}} gives a very brief overview of K-theory, Elliott's invariant, and the role of the UCT, especially for stable uniqueness results. Section {\ref{section-qd}} is devoted to quasidiagonality, and a rough outline of the main theorem of \cite{TWW}. In Section~{\ref{sectionRosenberg}} we revisit Rosenberg's conjecture on the connection between amenability and quasidiagonality. Section~{\ref{sectionTW}} summarises what is known and what is not known about the Toms--Winter conjecture for simple, unital, nuclear C$^*$-algebras. Finally, Section~{\ref{section-classification}} highlights the state of the art of Elliott's classification  programme.

I am indebted to Hannes Thiel, Aaron Tikuisis, and Stuart White for looking carefully at an earlier manuscript.

\section{Internal approximation: nuclearity and exactness}
\label{section1}

\bn
\label{def-cp}
C$^*$-algebras form a category, with the most natural choice for morphisms being $^*$-homomorphisms. (It follows from spectral theory that these are automatically norm-decreasing, hence continuous.) 

Another important class of morphisms consists of completely positive maps (we write c.p., or c.p.c.\ if they are also contractive). These are linear and $^*$-preserving, and send positive elements to positive elements, even after amplification with matrix algebras. By Stinespring's theorem, every c.p.\ map can be written as a compression of a $^*$-homomorphism. More precisely, a map $\varphi: A \longrightarrow B$ is completely positive if and only if $B$ embeds into another C$^*$-algebra $C$ and if there are a $^*$-homomorphism $\pi:A \longrightarrow C$ and some $h \in C$ such that 
\[
\varphi(\, .\,) = h^* \pi(\, .\, ) h.
\]   
\en

\bn
\label{def-order-zero}
A completely positive map $\varphi: A \longrightarrow B$ is order zero if it preserves orthogonality, i.e., whenever $a_1, \, a_2 \in A$ satify $a_1a_2 = 0$, then $\varphi(a_1) \varphi(a_2) = 0$. By \cite{WZ:MJM}, c.p.\ order zero maps are precisely those for which there is a Stinespring dilation such that $h \in C$ is positive and commutes with $\pi(A)$. 

As a consequence of this structure theorem, there is a bijective correspondence between c.p.c.\ order zero maps $A \longrightarrow B$ and $^*$-homomorphisms $\mathrm{C}_0((0,1], A) \longrightarrow B$. Moreover, one can use functional calculus on the commutative C$^*$-algebra generated by $h$ (which in fact is a quotient of $\mathrm{C}_0((0,1])$), to define a notion of functional calculus for c.p.c.\ order zero maps; cf.\ \cite{WZ:MJM}.
\en

\bn
\label{def-CPAP}
\begin{ndefn}
A $\mathrm{C}^*$-algebra $A$ has the completely positive approximation property, if the following holds: 

For any finite subset $\mathcal{F} \subset A$ and any tolerance $\epsilon>0$, there is a diagram
\begin{equation}
\label{approximation-diagram}
A \stackrel{\psi}{\longrightarrow} F \stackrel{\varphi}{\longrightarrow} A
\end{equation}
with $F$ a finite dimensional $\mathrm{C}^*$-algebra and completely positive contractive  maps $\psi$ and $\varphi$, such that $\varphi \psi$ agrees with the identity map up to $\epsilon$ on $\mathcal{F}$, in short
\[
\varphi \psi =_{\mathcal{F},\epsilon} \mathrm{id}_A ,
\mbox{ i.e., }
\|a - \varphi\psi(a)\| < \epsilon \mbox{ for all } a \in \mathcal{F}.
\]
\end{ndefn}

\noindent
We have asked the maps $\psi$ and $\varphi$ to be contractions. One could also ask them to be just bounded. As long as the norm bound is uniform, the resulting definitions will be equivalent. 
\en

\bn
A C$^*$-algebra $A$ is nuclear if, for every other C$^*$-algebra $B$, there is only one C$^*$-norm on the algebraic tensor product $A \odot B$; equivalently, the maximal and minimal tensor products of $A$ and $B$ agree. $A$ is exact if taking the minimal tensor product with another C$^*$-algebra $B$ is an exact functor for any $B$. Since the maximal tensor product has this property, nuclear C$^*$-algebras are automatically exact. 

Choi and Effros proved in \cite{CE:CPAP} (and Kirchberg in \cite{Kir:CPAP}) that a C$^*$-algebra is nuclear if and only if it has the completely positive approximation property. By the Choi--Effros lifting theorem, c.p.c.\ maps from nuclear C$^*$-algebras into quotient C$^*$-algebras always admit c.p.c.\ lifts. 

Kirchberg's $\mathcal{O}_2$ embedding theorem says that any separable exact C$^*$-algebra can be embedded into the Cuntz algebra $\mathcal{O}_2$. Separable nuclear C$^*$-algebras are precisely those which in addition are the images of conditional expectations on $\mathcal{O}_2$.   
\en

\bn 
\label{def-dimnuc}
In \cite{KirWin:dr}, Eberhard Kirchberg and I defined a notion of covering dimension for C$^*$-algebras which is based on \ref{def-CPAP} and uses order zero maps to model disjointness of open sets in a noncommutative situation. This notion, called decomposition rank, was generalised in \cite{WinZac:dimnuc} by Joachim Zacharias and myself. Here are the precise definitions.

\begin{ndefn}
A $\mathrm{C}^*$-algebra $A$ has nuclear dimension at most $d$, $\mathrm{dim}_{\mathrm{nuc}} \,A \le d$, if the following holds: 

For any finite subset $\mathcal{F} \subset A$ and any tolerance $\epsilon>0$, there is a diagram
\[
A \stackrel{\psi}{\longrightarrow} F \stackrel{\varphi}{\longrightarrow} A
\]
with $F$ a finite dimensional $\mathrm{C}^*$-algebra, $\psi$ completely positive contractive and $\varphi$ completely positive, such that 
\begin{itemize}
\item[(i)] $\|a - \varphi\psi(a)\| < \epsilon$ for all $a \in \mathcal{F}$,
\item[(ii)] there is a decomposition $F = F^{(0)} \oplus \ldots \oplus F^{(d)}$ such that each $\varphi^{(i)}:= \varphi|_{F^{(i)}}$ is c.p.c.\ order zero.
\end{itemize}
If, moreover, the maps $\varphi$ can be chosen to be contractive as well, then we say $A$ has decomposition rank at most $d$, $\mathrm{dr}\, A \le d$.
\end{ndefn}

Both of these concepts generalise covering dimension for locally compact spaces. The values are zero precisely for AF algebras, and they can be computed (or at least bounded from above) for many concrete examples. 

Note that for nuclear dimension, the maps $\varphi$ are sums of $d+1$ contractions, hence are uniformly bounded. Therefore, both finite nuclear dimension and finite decomposition rank imply nuclearity. The two concepts, as similar as they look, are genuinely different. In particular, unlike decomposition rank, nuclear dimension may be finite also for infinite C$^*$-algebras such as the Toeplitz algebra or the Cuntz algebras. The problem of characterising the difference turned out to be close to the heart of the subject. It was shown in \cite{KirWin:dr} and \cite{WinZac:dimnuc}, respectively, that the maps $\psi$ can be arranged to be approximately multiplicative for decomposition rank, and approximately order zero for nuclear dimension. In the former case, this shows that finite decomposition rank implies quasidiagonality; cf.\ Section~{\ref{section-qd}}. We will see that this is close to nailing down the difference between decomposition rank and nuclear dimension precisely.
\en

\bn
\label{convex-approximations}
The approximations of \ref{def-dimnuc} are more rigid than those of \ref{def-CPAP}. This means for some purposes they are more useful, but we also know that not all nuclear C$^*$-algebras have finite nuclear dimension. Building on \cite{HKW:Advances} (and extending \cite{CE:CPAP}), Brown--Carrion--White in \cite{BCW:Abel} gave a refined version of the completely positive approximation property, which asks the involved maps to be somewhat more rigid. In particular, the `downwards' maps $\psi$ can be taken to be approximately order zero, and the `upwards' maps $\varphi$ to be sums of honest order zero maps. The precise statement is as follows.

\begin{ntheorem}
A $\mathrm{C}^*$-algebra $A$ is nuclear if and only if the following holds: 

For any finite subset $\mathcal{F} \subset A$ and any tolerance $\epsilon>0$, there is a diagram
\[
A \stackrel{\psi}{\longrightarrow} F \stackrel{\varphi}{\longrightarrow} A
\]
with $F$ a finite dimensional $\mathrm{C}^*$-algebra and c.p.c.\  maps $\psi$ and $\varphi$, such that 
\begin{itemize}
\item[(i)] $\|a - \varphi\psi(a)\| < \epsilon$ for all $a \in \mathcal{F}$,
\item[(ii)] $\|\psi(a) \psi(b)\| < \epsilon$ whenever $a,b \in \mathcal{F}$ satisfy $ab = 0$,
\item[(iii)] there is a decomposition $F = F^{(0)} \oplus \ldots \oplus F^{(k)}$, such that the restrictions $\varphi|_{F^{(i)}}$ all have order zero, and such that $\sum_{i=1}^k \|\varphi|_{F^{(i)}}\| \le 1$.   
\end{itemize}
\end{ntheorem}

Compared to \ref{def-dimnuc}, in this statement the number of summands (of which I think as colours) is not uniformly bounded, but unlike in the original completely positive approximation property of \ref{def-CPAP} one still has individual order zero maps. As an extra bonus, one can arrange the norms to add up to one, or, upon normalising, one can think of the maps $\varphi$ as \emph{convex} combinations of c.p.c.\ order zero maps. This kind of approximation is a little subtle to write down explicitly, even when $A$ is a commutative C$^*$-algebra like $\mathrm{C}([0,1])$. (In particular, the number of colours in this setup will typically become very large; this is because in the proof at some point one has to pass from weak$^*$ to norm approximations via a convexity argument.) 
\en

\bn
One can think of the approximations of \ref{def-dimnuc} and \ref{convex-approximations} as \emph{internal} in the following sense:

By \ref{def-order-zero}, each order zero map corresponds to a $^*$-homomorphism from the cone over the domain C$^*$-algebra, which -- for each matrix summand -- essentially is given by an embedding of an algebra like $\mathrm{C}_0(X \setminus \{0\}) \otimes M_k$. Here, $X \subset [0,1]$ is a compact subset which is just the spectrum of the positive contraction $\varphi(1_{M_k}) \in A$. On the other hand, one may approximate $X \setminus\{0\}$ by the union of at most one half-open interval, finitely many  closed intervals, and finitely many points. Now one can use order zero functional calculus to slightly modify $\varphi|_{M_k}$ in such a way that the image sits in an honest subalgebra of $A$ which (after rescaling and relabeling the involved intervals) is isomorphic to 
\[
(\mathrm{C}_0((0,1]) \otimes M_k) \oplus (\mathrm{C}([0,1]) \otimes M_k) \oplus \ldots \oplus (\mathrm{C}([0,1]) \otimes M_k) \oplus M_k \oplus \ldots \oplus M_k.
\]
The overall $\varphi$ then maps into a (non-direct) sum of such subalgebras; with a little extra effort one can describe the map $\psi$ in terms of associated conditional expectations, and of course one can also keep track of the convex coefficients in \ref{convex-approximations}. It is usually more practical to write c.p.\ approximations like in \eqref{approximation-diagram}, but I often do find it useful to think of them as genuinely internal.
\en

\section{K-theory, the UCT, and stable uniqueness}
\label{section2}

\bn 
K-theory for C$^*$-algebras is a generalisation of topological K-theory; it is \emph{the} homology theory which is at the same time homotopy invariant, half-exact, and compatible with stabilisation. For a (say unital) C$^*$-algebra $A$, $\mathrm{K}_0(A)$ may be defined in terms of equivalence classes of projections in matrix algebras over $A$, or in terms of equivalence classes of finitely generated projective modules over $A$. (One first arrives at a semigroup, whose Grothendieck group is an ordered abelian group defined to be $\mathrm{K}_0$.) $\mathrm{K}_1(A)$ may be defined in terms of equivalence classes of unitaries, or by taking $\mathrm{K}_0$ of the suspension. Like for topological $\mathrm{K}$-theory, one has Bott periodicity, so $\mathrm{K}_{*+2} (A) \cong \mathrm{K}_{*} (A)$. 

Every tracial state on $A$ naturally induces a state (i.e., a positive real valued character) on the ordered $\mathrm{K}_0$-group. 
\en

\bn
For $A$ simple and unital, the Elliott invariant consists of the ordered $\mathrm{K}_0$-group (together with the position of the class of the unit), $\mathrm{K}_1$, and of the trace space together with the pairing with $\mathrm{K}_0$,
\[
\mathrm{Ell}(A) := (\mathrm{K}_0(A), \mathrm{K}_0(A)_+,[1_A]_0,\mathrm{K}_1(A), \mathrm{T}(A), \mathrm{r}_A: \mathrm{T}(A) \longrightarrow \mathrm{S}(\mathrm{K}_0(A))).
\]   
$\mathrm{Ell}(\,.\,)$ is a functor in a natural manner. 

We say a class $\mathcal{E}$ of simple unital C$^*$-algebras is classified by the Elliott invariant, if the following holds:

Whenever $A$, $B$ are in $\mathcal{E}$, and there is an isomorphism between $\mathrm{Ell}(A)$ and $\mathrm{Ell}(B)$, then there is an isomorphism between the algebras lifting the isomorphism of invariants. We will see in Section~{\ref{section-classification}} that this actually happens, and in great generality.
\en

\bn
Kasparov's $\mathrm{KK}$-theory is a bivariant functor from $\mathrm{C}^*$-algebras to abelian groups which is contravariant in the first and covariant in the second variable. It has similar abstract properties as $\mathrm{K}$-theory, and we have $\mathrm{K}_*(A) \cong \mathrm{KK}_*(\mathbb{C},A)$. 

Rosenberg and Schochet in \cite{RS:DMJ} studied the sequence
\begin{equation}
\label{UCT-sequence}
0 \longrightarrow \mathrm{Ext}^1(\mathrm{K}_*(A),\mathrm{K}_{*+1}(B)) \longrightarrow \mathrm{KK}(A,B) \longrightarrow \mathrm{Hom}(\mathrm{K}_*(A), \mathrm{K}_*(B)) \longrightarrow 0.
\end{equation}
A separable C$^*$-algebra $A$ is said to \emph{satisfy the universal coefficient theorem}  (UCT for short), if the sequence \eqref{UCT-sequence} is exact for every $\sigma$-unital C$^*$-algebra $B$. It follows from \cite{RS:DMJ} (cf.\ \cite[Theorem~23.10.5]{B:Kthy}), that $A$ satisfies the UCT precisely if it is $\mathrm{KK}$-equivalent to an abelian C$^*$-algebra.

The UCT problem asks whether \emph{all} separable \emph{nuclear} C$^*$-algebras satisfy the UCT. This is perhaps the most important open question about nuclear C$^*$-algebras.  
\en

\bn
\label{su-theorem}
The sequence \eqref{UCT-sequence} allows one to lift homomorphisms between $\mathrm{K}$-groups to $\mathrm{KK}$-elements. The latter are already a little closer to $^*$-homomorphisms between the C$^*$-algebras, but to get there one needs fairly precise control over the extent to which $^*$-homomorphisms are determined by their $\mathrm{KK}$-classes. This is often done by so-called stable uniqueness theorems, as developed in particular by Lin, Dadarlat and Eilers, and others. Let us state here a slightly simplified version of \cite[Theorem~4.5]{DE:PLMS}.  

\begin{ntheorem}
Let $A$, $B$ be unital $\mathrm{C}^*$-algebras with $A$ separable and nuclear. Let $\iota:A \longrightarrow B$ be a unital $^*$-homomorphism which is totally full, i.e., for every nonzero positive element of $A$, its image under $\iota$ generates all of $B$ as an ideal.
Let $\phi, \psi:A \longrightarrow B$ be unital $^*$-homomorphisms such that $\mathrm{KK}(\phi) = \mathrm{KK}(\psi)$.

Then, for every finite subset $\mathcal{G} \subset A$ and every $\delta > 0$ there are $n \in \mathbb{N}$ and a unitary $u \in M_{n+1}(B)$ such that
\[
\| u (\phi(a) \oplus \iota^{\oplus n}(a)) u^* - (\psi(a) \oplus \iota^{\oplus n}(a)) \| < \delta \mbox{ for all } a \in \mathcal{G}.
\]
\end{ntheorem}
\en

\bn
\label{su-theorem-UCT}
In the theorem above, the number $n$ depends on $\mathcal{G}$ and on $\delta$, but also on the maps $\phi$, $\psi$, and $\iota$. However, in applications one often cannot specify these maps beforehand. Dadarlat and Eilers in \cite{DE:PLMS} have found a way to deal with this issue (their original result only covers simple domains, but it can be pushed to the non-simple situation as well; cf.\ \cite[Lemma~5.9]{Lin:asu} or \cite[Theorem~3.5]{TWW}). The idea is it to assume that $n$ cannot be chosen independently of the maps, and then to construct sequences of maps  which exhibit this behaviour. Now regard these sequences as product maps, and apply the original Theorem~\ref{su-theorem} to arrive at a contradiction. To this end, it is important to keep control over the $\mathrm{KK}$-classes of the product maps --- which is not easy, since $\mathrm{KK}$-theory is not compatible with products in general. At this point the UCT saves the day, since (at least for the algebras involved) it guarantees that the map
\[
\textstyle
\mathrm{KK}\big(A,\prod_\mathbb{N} B_n\big) \longrightarrow \prod_\mathbb{N} \mathrm{KK}(A, B_n)
\]
is injective. Very roughly, if two sequences of KK-elements on the right hand side agree, they are connected by a sequence of homotopies. But since there is no uniform control over the length of these, it is not clear how to combine them to a single homotopy on the left hand side, at least not for general $A$. On the other hand, one can do this if the domain algebra $A$ has some additional geometric structure --- e.g., if it is commutative. But then of course it also suffices if $A$ is $\mathrm{KK}$-equivalent to a commutative C$^*$-algebra, i.e., if it satisfies the UCT. 
\en

\section{External approximation: quasidiagonality}
\label{section-qd}

\bn
Halmos defined a set $S \subset \mathcal{B}(\mathcal{H})$ of operators on a Hilbert space to be quasidiagonal if there is an increasing net of finite rank projections converging strongly to the identity operator on the Hilbert space, such that the projections approximately commute with elements of $S$. One then calls a C$^*$-algebra quasidiagonal if it has a faithful representation on some Hilbert space, such that the image forms a quasidiagonal set of operators in Halmos' sense. Voiculescu in \cite[Theorem~1]{Voi:Duke} rephrased this in a way highlighting quasidiagonality as an external approximation property. 

\begin{ntheorem}
A $\mathrm{C}^*$-algebra $A$ is quasidiagonal if, for every finite subset $\mathcal{F}$ of $A$ and $\epsilon > 0$, there are a matrix algebra $M_k$ and a c.p.c.\  map $\psi:A \longrightarrow M_k$ such that
\begin{itemize}
\item[(i)] $\|\psi(ab) - \psi(a) \psi(b)\| < \epsilon$ for all $a,b \in \mathcal{F}$,
\item[(ii)] $\|\psi(a)\| > \|a\| - \epsilon$ for all $a \in \mathcal{F}$.
\end{itemize}
\end{ntheorem}
\en

\bn
The maps of the theorem above may be thought of as approximate finite dimensional representations. This point of view has still not been fully exploited, partly because C$^*$-algebras are not so accessible to representation theoretic methods. It has, on the other hand, turned out to be fruitful to think of quasidiagonality as an embeddability property. Let $\mathcal{Q}$ denote the universal UHF algebra, $\mathcal{Q} = M_2 \otimes M_3 \otimes \ldots$, so that each matrix algebra embeds unitally into $\mathcal{Q}$. Then, a separable C$^*$-algebra $A$ is quasidiagonal if and only if there is a commuting diagram of the form
\[
\xymatrix{
&& \prod_\mathbb{N} \mathcal{Q} \ar[d] \\
A  \ar[rr]^{\bar{\psi}} \ar[urr]^{\tilde{\psi}}  && \prod_\mathbb{N} \mathcal{Q}/\bigoplus_\mathbb{N} \mathcal{Q} 
}
\]
with $\bar{\psi}$ an injective $^*$-homomorphism and $\tilde{\psi}$ a completely positive contraction. If, in addition, $A$ is nuclear, then the lift $\tilde{\psi}$ always exists by the Choi--Effros lifting theorem. Moreover, one may replace the sequence algebra $\prod_\mathbb{N} \mathcal{Q}/\bigoplus_\mathbb{N} \mathcal{Q}$ by an ultrapower $\mathcal{Q}_\omega$.
As a result, a separable nuclear C$^*$-algebra $A$ is quasidiagonal if and only if there is an embedding
\begin{equation}
\label{qd-Q-embedding}
\kappa: A \longrightarrow \mathcal{Q}_\omega.
\end{equation}
\en

\bn
\label{QDQ}
Every quasidiagonal C$^*$-algebra is stably finite, i.e., neither the algebra nor any of its matrix amplifications contains a projections which is Murray--von Neumann equivalent to a proper subprojection (this is a finiteness condition, reminiscent of the absence of paradoxical decompositions). The quasidiagonality question (QDQ) asks whether this is the only obstruction, at least in the nuclear case.

\begin{nquestion} 
(QDQ) Is every stably finite nuclear C$^*$-algebra quasidiagonal?
\end{nquestion}

After being around for some time this was first put in writing by Blackadar and Kirchberg in \cite{BlaKir:MathAnn}. There is a range of variations as discussed in \cite{Win:Abel}.
\en

\bn
By \cite{O:JMSUT}, like $\mathcal{Q}$ itself, the ultrapower $\mathcal{Q}_\omega$ also has a unique tracial state $\tau_{\mathcal{Q}_\omega}$. The composition $\tau_{\mathcal{Q}_\omega} \circ \kappa$ is a positive tracial functional on $A$. Whenever this is nonzero one may rearrange both $\kappa$ and its lift so that $\tau_{\mathcal{Q}_\omega} \circ \kappa$ is a tracial state on $A$. A tracial state $\tau$ which arises in this manner is called a quasidiagonal trace: 
\begin{equation*}
\label{qd-trace-diagram}
\xymatrix{
A  \ar[rr]^{\exists \kappa} \ar[drr]_{\tau}  && \mathcal{Q}_\omega \ar[d]^{\tau_{\mathcal{Q}_\omega}} \\
&& \mathbb{C}
}
\end{equation*}
It is common to drop the injectivity requirement on $\kappa$ in \eqref{qd-Q-embedding} in this context; this is largely for notational convenience since otherwise 
one would often have to factorise through the quotient by the trace kernel ideal.

A natural refinement of \ref{QDQ} is QDQ for traces; cf.\ \cite{BrOz,Win:Abel}. Just as QDQ, this has been around for a while; to the best of my knowledge it appeared in Nate Brown's \cite{B:MAMS} for the first time explicitly. It became a quite crucial topic for \cite{BBSTWW:arXiv} and \cite{TWW}, and also for Elliott's classification programme, as we will see below. 

\begin{nquestion}
Is every tracial state on a nuclear C$^*$-algebra quasidiagonal?
\end{nquestion}
\en

\bn
The fact that \emph{unital} quasidiagonal C$^*$-algebras always have at least one quasi\-diagonal trace was first observed by Voiculescu in \cite{V:IEOT}. On the other hand, an arbitrary embedding $\kappa: A \longrightarrow \mathcal{Q}_\omega$ may well end up in the trace kernel ideal of $\tau_{\mathcal{Q}_\omega}$, so that the composition $\tau_{\mathcal{Q}_\omega} \circ \kappa$ vanishes. For embeddings of cones as in \cite{Voi:Duke} this will always be the case. In that paper Voiculescu showed that quasidiagonality is homotopy invariant and concluded that cones and suspensions of arbitrary (say separable) C$^*$-algebras are always quasidiagonal. The method is completely general, but it does not allow to keep track of tracial states.

For cones over \emph{nuclear} C$^*$-algebras, one can say more: \cite[Lemma~2.6]{TWW} introduced a way of mapping a cone over the nuclear C$^*$-algebra $A$ to $\mathcal{Q}_\omega$ while at the same time controlling a prescribed trace on $A$. More precisely:

\begin{nlemma}
Let $A$ be a separable nuclear $\mathrm{C}^*$-algebra with a tracial state $\tau_A$. Then, 
\begin{enumerate}
\item[(i)]
there is a c.p.c.\ order zero map 
\[
\varphi:A \longrightarrow \mathcal{Q}_\omega
\]
such that $\tau_A = \tau_{\mathcal{Q}_\omega} \circ \varphi$,
and
\item[(ii)] there is a $^*$-homomorphism 
\[
\acute{\Lambda}: \mathrm{C}_0((0,1]) \otimes A \longrightarrow  \mathcal{Q}_\omega 
\]
such that $\tau_{\lambda} \otimes \tau_A = \tau_{\mathcal{Q}_\omega} \circ \acute{\Lambda}$, where $\tau_{\lambda}$ denotes the Lebesgue integral on $\mathrm{C}_0((0,1])$.
\end{enumerate}
\end{nlemma}

Let us have a quick glance at the proof. Since $\mathcal{Q}$ is tensorially self-absorbing and since it is not very hard to find an embedding $\lambda$ of the interval into $\mathcal{Q}$ in a Lebesgue trace preserving way, (ii) follows from (i) 
by extending the c.p.c.\ order zero map
\[
A \stackrel{\lambda(\mathrm{id}_{(0,1]}) \otimes \varphi}{\longrightarrow}  
\mathcal{Q} \otimes \mathcal{Q}_\omega \longrightarrow (\mathcal{Q} \otimes \mathcal{Q})_\omega \cong \mathcal{Q}_\omega
\]
to a $^*$-homomorphism defined on the cone,
\[
\acute{\Lambda}: \mathrm{C}_0((0,1]) \otimes A \longrightarrow  \mathcal{Q}_\omega. 
\]

For (i), for the sake of simplicity let us assume that $\tau_A$ is extremal. Then since $A$ is nuclear, by Connes' work \cite{C:Ann}, the weak closure of $A$ in the GNS representation $\pi_{\tau_A}$ for $\tau_A$ is the hyperfinite $\mathrm{II}_1$ factor $\mathcal{R}$. It follows from the Kaplansky density theorem that there is a surjection from $\mathcal{Q}_\omega$ onto $\mathcal{R}_\omega$. Now again by nuclearity, the Choi--Effros lifting theorem yields a c.p.c.\ lift $\widetilde{\varphi}$ of $\pi_{\tau_A}$:
\[
\xymatrix{
&&&  \mathcal{Q}_\omega \ar[d]^{\mathrm{q}} \\
A  \ar[rr]_{\pi_{\tau_A}} \ar[urrr]^{\widetilde{\varphi}}&& \mathcal{R} \ar@{^{(}->}[r]    & \mathcal{R}_\omega 
}
\]
This lift $\widetilde{\varphi}$ has no reason to be order zero. However, for any approximate unit $(e_\lambda)_\Lambda$ of the kernel of the quotient map $\mathrm{q}$, the maps $\widetilde{\varphi}_\lambda := (1-e_\lambda)^{1/2}\widetilde{\varphi}(\, . \,)(1-e_\lambda)^{1/2}$ will lift $\pi_{\tau_A}$ as well --- and if one takes the approximate unit to be quasicentral with respect to $Q_\omega$ those maps are at least \emph{approximately} order zero. Now use separability of $A$ and a `diagonal sequence argument' to turn the $\widetilde{\varphi}_\lambda$ into an honest order zero lift $\varphi$. This type of diagonal sequence argument appears inevitably when working with sequence algebras. In this case one can run it more or less by hand, but a better, and more versatile way to implement it in a C$^*$-algebra context is Kirchberg's $\epsilon$-test; cf.\ \cite[Lemma~A1]{K:Abel}.
\en

\bn
Let us take another look at the lemma above when $A$ is unital. In this case, we have an embedding 
\[
\acute{\lambda}:= \acute{\Lambda}|_{\mathrm{C}_0((0,1]) \otimes 1_A}: \mathrm{C}_0((0,1]) \longrightarrow \mathcal{Q}_\omega
\]
of the cone into  $\mathcal{Q}_\omega$. One may unitise this map to arrive at an embedding
\[
\bar{\lambda}: \mathrm{C}([0,1]) \longrightarrow \mathcal{Q}_\omega
\] 
which still induces the Lebesgue integral when composed with $\tau_{\mathcal{Q}_\omega}$. If only we could extend this map $\bar{\lambda}$ to $\mathrm{C}([0,1]) \otimes A$, then  this would immediately prove quasidiagonality of the trace $\tau_A$. Of course such an extension seems far too much to ask for, but it is not completely unreasonable either: The map $\bar{\lambda}$ restricts to the embeddings 
\[
\acute{\lambda}: \mathrm{C}_0((0,1]) \longrightarrow \mathcal{Q}_\omega \; \mbox{ and } \;\grave{\lambda}: \mathrm{C}_0([0,1)) \longrightarrow \mathcal{Q}_\omega.
 \]
 Now since the Lebesgue integral is symmetric under flipping the interval, we see that it agrees with both maps $\tau_{\mathcal{Q}_\omega} \circ \acute{\lambda}$ and $\tau_{\mathcal{Q}_\omega} \circ \grave{\lambda} \circ \mathrm{flip}$ (where flip denotes the canonical isomorphism between $\mathrm{C}_0((0,1])$ and $\mathrm{C}_0([0,1))$. Moreover, by \cite{CE:IMRN} this is enough to make the maps  $\acute{\lambda}\circ \mathrm{flip}$ and $\grave{\lambda}$ approximately unitarily equivalent --- and again by a diagonal sequence argument one can even make them honestly unitarily equivalent, i.e., one can find a unitary $u \in \mathcal{Q}_\omega$ such that
 \begin{equation}
 \label{u-flip}
 \grave{\lambda}  (\, .\,) = u\, \acute{\lambda}(\mathrm{flip}(\,.\,)) \,u^*.
 \end{equation}
Now this map can clearly be extended to all of $\mathrm{C}_0([0,1)) \otimes A$ by setting
\[
\grave{\Lambda}:= u\, \acute{\Lambda}((\mathrm{flip} \otimes \mathrm{id}_A)(\,.\,)) \,u^*.
\] 
At this point we have two maps 
\[
\acute{\Lambda}: \mathrm{C}_0((0,1]) \otimes A \longrightarrow \mathcal{Q}_\omega \; \mbox{ and } \;\grave{\Lambda}: \mathrm{C}_0([0,1)) \otimes A \longrightarrow \mathcal{Q}_\omega,
\]
which we would like to `superpose' to a map defined on  $\mathrm{C}([0,1]) \otimes A$. This can be done by means of a $2\times2$ matrix trick, involving the unitary $u$ and rotation `along the interval'. The result will be a c.p.c.\ map 
\[
\bar{\Lambda}: \mathrm{C}([0,1]) \otimes A \longrightarrow M_2(\mathcal{Q}_\omega)
\]
which will map $1_{[0,1]}$ to a projection of trace $1/2$. However, to arrive at quasidiagonality, $\bar{\Lambda}$ would also have to be multiplicative. This will indeed happen provided one can in addition choose the unitary $u$ to implement the flip on all of the suspension $\mathrm{C}_0((0,1)) \otimes A$, or equivalently, to satisfy \eqref{u-flip} as well as   
 \begin{equation}
 \label{flip-grave-acute}
\grave{\Lambda}(f \otimes a) = u\, \acute{\Lambda}(f  \otimes a) \,u^* = \acute{\Lambda}(f  \otimes a)
\end{equation} 
for all $a \in A$ and for all $f \in \mathrm{C}_0((0,1))$ which are symmetric in the sense that $\mathrm{flip}(f) = f$. In other words, we have to implement the flip on $\mathrm{C}_0((0,1))$ in the relative commutant of a certain suspension over $A$.

It is a lot to ask for such a unitary to begin with, and even approximate versions are just as  hard to achieve, since in an ultrapower approximately implementing \eqref{u-flip} and \eqref{flip-grave-acute} will be as good as implementing them exactly. On the other hand, the domains of our maps are cones, or suspensions embedded in cones, hence zero-homotopic, and so there is no obstruction in $\mathrm{K}$-theory to finding such a $u$. Luckily, there are powerful techniques from C$^*$-algebra classification in place which do allow to compare maps when they sufficiently agree on $\mathrm{K}$-theory; cf.\ \cite{Lin:stableJOT, DE:PLMS}, or \ref{su-theorem} above. These require the target algebra to be `admissible' (which is the case here), but there is also a price to pay: with these stable uniqueness theorems, one can only compare maps up to (approximate) unitary equivalence after adding a `large' map to both sides. This largeness can be measured numerically using the trace of the target algebra (which recovers the trace $\tau_A$ on $A$ via $\Lambda$; at this point it is important that $\tau_A$ is faithful). In general the largeness constant depends on the algebras involved, but also on the maps. In \cite{TWW} we found a way to use large multiples of the original maps as correcting summands. It is then important for the largeness constant to not depend on the maps involved. As pointed out in \ref{su-theorem-UCT}, Dadarlat and Eilers have indeed developed such a stable uniqueness theorem which works when the domain (in our case the suspension $\mathrm{C}_0((0,1)) \otimes A$) in addition satisfies the UCT. 
 
 This is all made precise in \cite{TWW}, which also contains an extensive sketch of the proof (a slightly more informal sketch can be found in \cite{Win:Abel}). Here is the result.
 \en
 
 \bn
 \begin{ntheorem}
 \label{TWW-theorem}
 Every faithful trace on a nuclear $\mathrm{C}^*$-algebra satisfying the UCT is quasidiagonal.
\end{ntheorem}
 
In particular this answers the quasidiagonality question QDQ for UCT C$^*$-algebras with faithful traces. We will see some more consequences in the subsequent sections.
\en

\bn
In \cite{Gab:qd-exact}, Gabe generalised the theorem above to the situation where $A$ is only assumed to be exact (but still satisfying the UCT), and the trace is amenable.

In \cite{Schaf:TWW}, Schafhauser gave a different, and shorter, proof, which replaces the stable uniqueness theorem of \cite{DE:PLMS} by a result from Elliott and Kucerovsky's \cite{EllKuc:PJM}. 
\en

\section{Rosenberg's conjecture: amenability}
\label{sectionRosenberg}

\bn
In the appendix of \cite{Hadwin:JOT}, Rosenberg observed that for reduced group C$^*$-algebras amenability and quasidiagonality are closely related.

\begin{nprop}
Let $G$ be a countable discrete group and suppose the reduced group $\mathrm{C}^*$-algebra $\mathrm{C}^*_{\mathrm{r}}(G) \subset \mathcal{B}(\ell^2(G))$ is quasidiagonal.
Then, $G$ is amenable.
\end{nprop}

A proof is not hard, and worth looking at (the one below can be extracted from \cite[Corollary~7.1.17, via Theorem~6.2.7]{BrOz}): For $g \in G$ let $u_g$ denote the image of $g$ in $\mathcal{B}(\ell^2(G))$ under the left regular representation. Let $(p_n)_\mathbb{N} \subset  \mathcal{B}(\ell^2(G))$ be a sequence of finite rank projections strongly converging to the identity and approximately commuting with $\mathrm{C}^*_{\mathrm{r}}(G)$. Compression with the $p_n$ yields unital c.p.\ maps 
\[
\varphi_n: \mathrm{C}^*_{\mathrm{r}}(G) \longrightarrow p_n \mathcal{B}(\ell^2(G)) p_n \cong M_{r_n}
\]
(where $r_n$ is just the rank of $p_n$). Since the $p_n$ approximately commute with elements of $\mathrm{C}^*_{\mathrm{r}}(G)$, these maps are approximately multiplicative, so that the limit map
\[
\textstyle
\varphi_\infty: \mathrm{C}^*_{\mathrm{r}}(G) \longrightarrow \prod_\mathbb{N} M_{r_n}/ \bigoplus_\mathbb{N} M_{r_n}
\]
is a $^*$-homomorphism. Note that each $\varphi_n$ extends to a unital c.p.\ map
\[
\bar{\varphi}_n:  \mathcal{B}(\ell^2(G)) \longrightarrow M_{r_n}
\]
by Arveson's extension theorem. Now by Stinespring's theorem, $\mathrm{C}^*_{\mathrm{r}}(G)$ sits in the multiplicative domain of the limit map 
\[
\textstyle
\bar{\varphi}_\infty:  \mathcal{B}(\ell^2(G)) \longrightarrow \prod_\mathbb{N} M_{r_n}/ \bigoplus_\mathbb{N} M_{r_n},
\]
which in particular means that for every $g \in G$ and every $x \in \mathcal{B}(\ell^2(G))$,
\[
\bar{\varphi}_\infty(u_g x) = \bar{\varphi}_\infty(u_g) \bar{\varphi}_\infty(x).
\]
Upon choosing a free ultrafilter $\omega \in \beta\mathbb{N} \setminus \mathbb{N}$, the canonical tracial states on the $M_{r_n}$, evaluated along $\omega$, yield a tracial state $\tau_\omega$ on $\prod_\mathbb{N} M_{r_n}/ \bigoplus_\mathbb{N} M_{r_n}$. Now for $x \in \mathcal{B}(\ell^2(G))$ we have
\begin{eqnarray*}
\tau_\omega \circ \bar{\varphi}_\infty(u_g x u_g^*) & = & \tau_\omega (\bar{\varphi}_\infty(u_g) \bar{\varphi}_\infty(x) \bar{\varphi}_\infty(u_g^*)) \\
& = & \tau_\omega ( \bar{\varphi}_\infty(u_g^*) \bar{\varphi}_\infty(u_g) \bar{\varphi}_\infty(x)) \\
& = & \tau_\omega ( \bar{\varphi}_\infty(u_g^*u_g) \bar{\varphi}_\infty(x)) \\
& = & \tau_\omega \circ \bar{\varphi}_\infty(x).
\end{eqnarray*}
This in particular holds for $x \in \ell^\infty(G)$ (regarded as multiplication operator), and we see that $\tau_\omega \circ \bar{\varphi}_\infty$ is a translation invariant state on $\ell^\infty(G)$. The existence of such an invariant mean is equivalent to $G$ being amenable.
\en

\bn
In the argument above, the $p_n$ approximately commute with elements of $\mathrm{C}^*_{\mathrm{r}}(G)$ \emph{in norm}. However, the construction almost forgets about this and really only requires the $p_n$ to approximately commute with $\mathrm{C}^*_{\mathrm{r}}(G)$ \emph{in trace}. The point is that $\tau_\omega \circ \bar{\varphi}_\infty$ is an \emph{amenable trace}, which is enough to show amenability of $G$; cf.\ \cite[Proposition~6.3.3]{BrOz}. If one conversely starts with an amenable group $G$ with a sequence of F{\o}lner sets $F_n$, and chooses $p_n \in \mathcal{B}(\ell^2(G))$ to be the associated finite rank projections, then the same construction as above will again yield an invariant mean. Of course the $p_n$ will approximately commute with $\mathrm{C}^*_{\mathrm{r}}(G)$ only in trace, and not in norm.

On the other hand, for $G = \mathbb{Z}$, one has $\mathrm{C}(S^1) \cong \mathrm{C}^*_{\mathrm{r}}(G) \subset \mathcal{B}(\ell^2(G))$, and it is well known that commutative C$^*$-algebras are quasidiagonal, not just in the abstract sense, but also when they are concretely represented on a Hilbert space. In this situation, one can even construct the quasi-diagonalising projections fairly explicity, from F{\o}lner sets, i.e., one can modify such F{\o}lner projections to make them approximately commute with $\mathrm{C}^*_{\mathrm{r}}(G)$ even in norm. (The idea is to `connect' the left and right hand sides of F{\o}lner sets inside the matrix algebra hereditarily generated by the projection.) 
\en

\bn 
The question of when (and how) one can find projections quasi-diagonalising $\mathrm{C}^*_{\mathrm{r}}(G)$ has turned out to be a hard one. Despite having ever so little evidence at hand at the time, Rosenberg did conjecture that amenable discrete groups are always quasidiagonal. He did not put the conjecture in writing in \cite{Hadwin:JOT}, but did promote the problem subsequently; see \cite{B:MAMS} and \cite{BrOz} for a more detailed discussion. 

The conjecture received attention by a number of researchers, and was indeed confirmed for larger and larger classes of amenable groups. These arguments often start with the abelian case and use some kind of bootstrap argument to reach more general classes of groups. The problem usually is that quasidiagonality does not pass to extensions.
\en

\bn
In \cite{ORS:GAFA}, Ozawa, R{\o}rdam and Sato proved Rosenberg's conjecture for elementary amenable groups. The latter form a bootstrap class, containing many but not all amenable groups (Grigorchuk's examples with exponential growth are amenable but not elementary amenable --- but their group C$^*$-algebras were already known to be quasidiagonal for other reasons). The argument of \cite{ORS:GAFA}  relies on methods and results from the classification programme for simple nuclear C$^*$-algebras. Therefore, already \cite{ORS:GAFA} factorises through a stable uniqueness result like \ref{su-theorem}. 
\en

\bn
Eventually, Rosenberg's conjecture was confirmed in full generality as a consequence of the main result from \cite{TWW}. 

\begin{ncor}
If $G$ is a discrete amenable group, then $\mathrm{C}^*_{\mathrm{r}}(G)$ is quasidiagonal.
\end{ncor}

First, it is well known (and not too hard to show) that the canonical trace on $\mathrm{C}^*_{\mathrm{r}}(G)$ is faithful. Then one consults Tu's work \cite{Tu:KT} to conclude that amenable group C$^*$-algebras always satisfy the UCT. Theorem~\ref{TWW-theorem} now says that $\mathrm{C}^*_{\mathrm{r}}(G)$ has a faithful quasidiagonal representation, i.e., it is quasidiagonal as an abstract C$^*$-algebra. But even in its concrete representation on $\ell^2(G)$ it is quasidiagonal; cf.\ \cite[Theorem~7.2.5]{BrOz}.\footnote{Note that countability / separability is not an issue since all properties involved can be tested locally.}
\en

\bn
The Corollary above does indeed settle Rosenberg's conjecture, but of course only in a very abstract manner. In particular, at this point there seems no way to exhibit quasi-diagonalising projections explicitly, starting, say, with a F{\o}lner system for $G$. 
\en

\section{Toms--Winter regularity}
\label{sectionTW}

\bn
Kirchberg's $\mathcal{O}_\infty$-absorption theorem says that a separable, simple, nuclear  C$^*$-algebra is purely infinite precisely if it absorbs the Cuntz algebra $\mathcal{O}_\infty$, $A \cong A \otimes \mathcal{O}_\infty$; see \cite{Kir:ICM}. Next to his $\mathcal{O}_2$-embedding theorem, this is one of the cornerstones for Kirchberg--Phillips classification; cf.\ \cite{R:Book} for an overview. 

At the time it was not at all clear whether one should expect a similar statement for stably finite C$^*$-algebras. We now know that the Jiang--Su algebra $\mathcal{Z}$ of \cite{JS:AJM} really is the right analogue of $\mathcal{O}_\infty$ in this context. Moreover, we know that pure infiniteness can be interpreted as a regularity property of the Cuntz semigroup (almost unperforation, to be more specific) in the absence of traces; cf.\ \cite{Ror:ICM}. On the other hand, the state of Elliott's classification programme in the early 2000s suggested that dimension type conditions should also play a crucial role. 
\en

\bn
\label{TW-conjecture}
In \cite{TomsWinter:VI}, Andrew Toms and I exhibited a class of inductive limit C$^*$-algebras for which  finite decomposition rank, $\mathcal{Z}$-stability, and almost unperforation of the Cuntz semigroup occur or fail simultaneously. This class (Villadsen algebras of the first type) was somewhat artificial, and a bit thin, but still large enough to prompt our conjecture that these three conditions should be equivalent for separable, simple, unital, nuclear and stably finite C$^*$-algebras. Once nuclear dimension was invented and tested, it became soon clear that the conjecture should be generalised to comprise both nuclear dimension and decomposition rank. The full version reads as follows.

\begin{nconj}
For a separable, simple, unital, nuclear $\mathrm{C}^*$-algebra $A \neq M_r$ the following are equivalent:
\begin{itemize}
\item[(i)] $A$ has finite nuclear dimension.
\item[(ii)] $A$ is $\mathcal{Z}$-stable.
\item[(iii)] $A$ has strict comparison of positive elements.
\end{itemize}
Under the additional assumption that $A$ is stably finite, one may replace {\rm (i)} by
\begin{itemize}
\item[(i')] $A$ has finite decomposition rank.
\end{itemize}
\end{nconj}
\en

\bn
\label{TW-whats-known}
I stated the above as a conjecture (perhaps in part for sentimental reasons), but after hard work by many people it is now almost a theorem (i.e., most of the implications have been proven in full generality). Let us recap what's known at this point.  

When $A$ has no trace, then (ii) $\Longleftrightarrow$ (iii) follows from Kirchberg's $\mathcal{O}_\infty$-absorption theorem as soon as one knows that an infinite exact C$^*$-algebra is $\mathcal{Z}$-stable if and only if it is $\mathcal{O}_\infty$-stable; see \cite{R:IJM}. 

In the finite case, (ii) $\Longrightarrow$ (iii) was shown by R{\o}rdam in \cite{R:IJM}. 

I showed first (i') $\Longrightarrow$ (ii) and then (i) $\Longrightarrow$ (ii) in \cite{W:Invent1} and \cite{W:Invent2}, respectively. 

In \cite{MS:Acta}, Matui and Sato showed (iii) $\Longrightarrow$ (ii) when $A$ has only one tracial state.  In each of \cite{S:Preprint2}, \cite{KR:Crelle}, and \cite{TWW:IMRN} this was generalised to the case where the tracial state space of $A$ (always a Choquet simplex) has compact and finite dimensional extreme boundary. 

In \cite{MS:DMJ}, Matui and Sato showed (ii) $\Longrightarrow$ (i') when $A$ has only finitely many extremal tracial states, and under the additional assumption that $A$ is quasidiagonal. In \cite{SWW:Invent}, (ii) $\Longrightarrow$ (i) was implemented in the case of a unique tracial state. In the six author paper \cite{BBSTWW:arXiv}, Joan Bosa, Nate Brown, Yasuhiko Sato, Aaron Tikuisis, Stuart White and myself showed (ii) $\Longrightarrow$ (i) when the tracial state space of $A$ has compact extreme boundary; (ii) $\Longrightarrow$ (i') was shown assuming in addition that every trace is quasidiagonal (by Theorem~\ref{TWW-theorem} this is automatic when $A$ satisfies the UCT). For these last results, the ground was prepared by Ozawa's theory of von Neumann bundles from \cite{O:JMSUT}. In upcoming work, Jorge Castillejos, Sam Evington, Aaron Tikuisis, Stuart White and I will show (ii) $\Longrightarrow$ (i) in full generality, and (ii) $\Longrightarrow$ (i') provided that every trace is quasidiagonal.  
\en

\bn
With all these results in place now, to sum up it is shorter to state what's not yet known: All we are missing is (iii) $\Longrightarrow$ (ii) for arbitrary trace spaces, and (i) $\Longrightarrow$ (i') without any quasidiagonality assumption. 

I am still quite optimistic about the latter statement. For the former one, I also remain positive, but every once in a while I'm tempted to travel back in time to replace condition (iii) by
\begin{itemize}
\item[(iii')] $A$ has strict comparison and has almost divisible Cuntz semigroup.
\end{itemize}
(This condition is equivalent to saying that the Cuntz semigroups of $A$ and $A \otimes \mathcal{Z}$ agree.) On the other hand, Thiel has recently shown that almost divisibility follows from strict comparison in the stable rank one case; see \cite{Thi:ranks}. Without this assumption, I wouldn't be too surprised if lack of divisibility was a new source of high dimensional examples in the spirit of \cite{T:Ann}.
\en

\bn
Ever since its appearance, Connes' classification of injective II$_1$ factors was an inspiration for the classification and structure theory of simple nuclear C$^*$-algebras. Once Conjecture~\ref{TW-conjecture} was formulated, it did not take that long to realise the surprising analogy with Connes' work, in particular \cite[Theorem~5.1]{C:Ann}. 
Roughly speaking, nuclearity on the C$^*$-algebra side corresponds to injectivity for von Neumann algebras, finite nuclear dimension to hyperfiniteness, $\mathcal{Z}$-stability to $\mathcal{R}$-absorption, i.e., being McDuff (cf.\ \cite{M:PLMS}), and strict comparison corresponds to comparison of projections. 

Matui and Sato in \cite{MS:Acta} and \cite{MS:DMJ} have taken this analogy to another level, by turning it into actual theorems. This trend was further pursued in \cite{SWW:Invent} and \cite{BBSTWW:arXiv}. I find these extremely convincing; by now I am even optimistic that eventually we will be able to view Connes' result and Conjecture~\ref{TW-conjecture} as incarnations of the same abstract theorem. 
\en

\section{Elliott's programme: classification}
\label{section-classification}

\bn
The first general classification result for nuclear C$^*$-algebras was probably Glimm's classification of UHF algebras in terms of supernatural numbers. Bratteli observed that one can do essentially the same for AF algebras using Bratteli diagrams, but it was Elliott who classified AF algebras in terms of their ordered $\mathrm{K}_0$-groups; \cite{E:JA}. More classification results for larger classes were picked up in the 1980s and early 1990s; these prompted Elliott to conjecture that separable, simple, nuclear C$^*$-algebra might be classifiable by K-theoretic invariants. (The precise form of the invariant underwent some adjustments as the theory and understanding of examples progressed.) Up to that point, all available results covered certain types of inductive limit C$^*$-algebras. Then Kirchberg opened the door to classification in a much more abstract context; \cite{Kir:ICM}. 
\en

\bn
I have already mentioned that Voiculescu showed quasidiagonality to be a homotopy invariant property. This in particular means that cones over separable C$^*$-algebras are quasidiagonal, because the former are contractible; since quasidiagonality passes to subalgebras, suspensions are quasidiagonal as well. Kirchberg used this statement to prove his celebrated $\mathcal{O}_2$-embedding theorem (cf.\ \cite{Kir:ICM}; see also \cite{R:Book}), which was a cornerstone for Kirchberg--Phillips classification of separable, nuclear, simple, purely infinite C$^*$-algebras. 

This was perhaps the earliest indication that quasidiagonality should be relevant for the classification of nuclear C$^*$-algebras, but of course in this situation we have quasidiagonal cones over C$^*$-algebras which are themselves very far from being quasidiagonal.   
\en

\bn
Soon after, however, Popa in \cite{Popa:PJM} carried over his local quantisation technique from von Neumann factors to simple C$^*$-algebras with traces and sufficiently many projections. The theorem roughly says that such C$^*$-algebras can be approximated locally by finite dimensional C$^*$-subalgebras. It is safe to say that this result kicked off the systematic use of quasidiagonality in the classification of \emph{stably finite} simple nuclear C$^*$-algebras. 
\en

\bn
\label{TAF-def}
Next, Lin modified Popa's local approximation to the effect that the approximating subalgebras are moreover required to be large in a certain sense; this can often be measured tracially, hence the name TAF (tracially approximately finite dimensional) C$^*$-algebras. Here is the precise definition.

\begin{ndefn}
A simple, unital $\mathrm{C}^*$-algebra $A$ is TAF, if the following holds: For every finite subset $\mathcal{F} \subset A$, $\epsilon>0$, and positive contraction $0 \neq e \in A$, there are a finite dimensional $\mathrm{C}^*$-subalgebra $F \subset A$ and a partial isometry $s \in A$ such that
\begin{itemize}
\item[(i)] $\mathrm{dist}(1_F a 1_F, F) < \epsilon$ for all $a \in \mathcal{F}$,
\item[(ii)] $\|1_F a - a 1_F\| < \epsilon$  for all $a \in \mathcal{F}$,
\item[(iii)] $s^*s = 1_A - 1_F$ and $ss^* \in \overline{eAe}$.
\end{itemize}
\end{ndefn}

The element $e$ dominates the complement of $1_F$, and therefore controls the size of $F$ inside $A$. In many situations, this can be done in terms of tracial states on $A$, hence the name \emph{tracially} AF.
\en

\bn
\label{TAF-theorem}
In \cite{Lin:Duke}, Lin managed to classify nuclear TAF algebras satisfying the UCT. The proof is inspired by Kirchberg--Phillips classification. 

\begin{ntheorem}
The class of all separable, simple, unital, nuclear, infinite dimensional, TAF, UCT $\mathrm{C}^*$-algebras is classified by the Elliott invariant.
\end{ntheorem}
\en

\bn
The theorem covers a fairly large class of stably finite C$^*$-algebras which is characterised abstractly (as opposed to the various classes of inductive limit type algebras handled earlier). The UCT hypothesis remains mysterious, but in applications, for example to transformation group C$^*$-algebras, this is often no issue since one can confirm the UCT directly. The scope of the theorem is nonetheless limited by the TAF assumption, which in particular requires the existence of many projections, and also the ordered $\mathrm{K}_0$-group to be weakly unperforated. Weakly unperforated $\mathrm{K}$-theory is implied by $\mathcal{Z}$-stability of the algebra; in \cite{MS:DMJ} it was shown that nuclear TAF algebras are indeed $\mathcal{Z}$-stable. In \cite{Win:localizing} it was shown that classification of $\mathcal{Z}$-stable C$^*$-algebras can be derived from classification of UHF-stable ones. This paved the road to applying Lin's TAF classification also in situations when the algebras contain no or only few projections. 

In \cite{GLN:arXiv}, Gong, Lin and Niu generalised Definition ~\ref{TAF-def} by admitting more general building blocks  than just finite dimensional C$^*$-algebras. At the same time, they managed to prove a classification result like Theorem~\ref{TAF-theorem} also in this context. This is a spectacular outcome, since the class covered by the result is no longer subject to $\mathrm{K}$-theoretic restrictions other than those implied by $\mathcal{Z}$-stability anyway. 

In \cite{EGLN:arXiv}, finally, it was shown that the UCT, together with finite nuclear dimension and quasidiagonality of all traces suffices to arrive at classification. 
In conjunction with \ref{TW-whats-known} and Theorem~\ref{TWW-theorem}, this confirms Elliott's conjecture in the $\mathcal{Z}$-stable, finite, UCT case. The $\mathcal{Z}$-stable, infinite, UCT case is precisely Kirchberg--Phillips classification. Moreover, R{\o}rdam's and Toms' examples in \cite{R:Acta} and in \cite{T:Ann} (inspired by Villadsen's \cite{V:JAMS}) have shown that $\mathcal{Z}$-stability cannot be disposed of for classification via the Elliott invariant. 

The upshot of this discussion is a classification result which---modulo the UCT problem---is as complete and final as can be. It is the culmination of decades of work, by many many hands. It reads as follows.

\begin{ntheorem}
The class of all separable, simple, unital, nuclear, $\mathcal{Z}$-stable UCT $\mathrm{C}^*$-algebras is classified by the Elliott invariant.
\end{ntheorem} 
\en

\nocite{EllToms:BAMS}
\nocite{Bro:symbiosis}
\nocite{Popa:JOT}
\nocite{Cun:CMP}
\nocite{EE:Ann}
\nocite{EGL:Invent}
\nocite{H:JFA}
\nocite{Haa:quasitraces}
\nocite{HigKas:Invent}

\bibliographystyle{amsplain}
\bibliography{ww}

\end{document}